\documentclass[11pt,a4paper]{article}

\usepackage{authblk}

\usepackage[T1]{fontenc}    
\usepackage{enumitem}
\usepackage{mathrsfs}
\usepackage{amsfonts}
\usepackage{amsthm}
\usepackage{authblk}
\usepackage{cite}
\usepackage{amsmath,amssymb,amscd,amsthm}
\usepackage{graphics}
\input xy
\xyoption{all}

\numberwithin{equation}{section}

\usepackage{scalefnt}   

\usepackage{comment}
\usepackage{setspace}

\usepackage[paperwidth=210mm,
paperheight=297mm,
a4paper,
twoside,
top=0.8in, left=0.9in, right=0.9in, bottom=1in]{geometry}

\usepackage[active]{srcltx}

\theoremstyle{definition}
\newtheorem{lemma}{Lemma}[section]

\newtheorem{prop}[lemma]{Proposition}

\newtheorem{defi}[lemma]{Definition}
\newtheorem{thm}[lemma]{Theorem}

\newtheorem{rmk}[lemma]{Remark}

\newcommand{\C}{\mathbb{C}}
\newcommand{\N}{\mathbb{N}}
\newcommand{\Z}{\mathbb{Z}}

\def\a{\alpha}

\usepackage{authblk} 
\usepackage{hyperref}
\usepackage{float}
\restylefloat{figure}

\allowdisplaybreaks

\begin{document}
	\title{Simple restricted modules of non-zero level over a deformed Heisenberg-Virasoro algebra}

    \author[1]{Shun Liu}
    
    \author[2]{Dashu Xu\thanks{Corresponding Author: dox@mail.ustc.edu.cn}}

    \affil[1]{School of Mathematical Sciences, Xiamen University, Xiamen 361005, Fujian, P. R. China}
    \affil[2]{School of Mathematical Sciences, University of Science and Technology of China, Hefei 230026, Anhui, P. R. China}

	\date{}
	\maketitle
	\begin{abstract}
    We study representations of a deformed Heisenberg-Virasoro algebra that does not admit a triangular decomposition. Despite this, its $\Z$-gradation allows the classification of simple restricted modules. We show that all such modules of non-zero level arise via induction from simple modules of finite-dimensional solvable Lie algebras.
    
	\bigskip
		
    \bigskip
    
		
	\end{abstract}
	
	\section{Introduction}
     Let $\mathbb{C}\left[t^{\pm1}, \frac{\mathrm{d}}{\mathrm{d}t}\right]$ denote the (associative) Weyl algebra. It contains a Lie subalgebra
    \[
    \mathcal{D} := \left\{ f(t)\frac{\mathrm{d}}{\mathrm{d}t} + g(t) \;\Big|\; f(t), g(t) \in \mathbb{C}[t^{\pm1}] \right\},
    \]
equipped with the Lie bracket
\[
\left[ f(t)\frac{\mathrm{d}}{\mathrm{d}t} + g(t), \, p(t)\frac{\mathrm{d}}{\mathrm{d}t} + q(t) \right] = \big( f(t)p'(t) - p(t)f'(t) \big)\frac{\mathrm{d}}{\mathrm{d}t} + \big( f(t)q'(t) - p(t)g'(t) \big).
\]
The Lie algebra $\mathcal{D}$ is referred to as the Lie algebra of differential operators on $\mathbb{C}[t^{\pm1}]$ of degree at most one\cite{M}.
     
     Denote the universal central extension of $\mathcal{D}$ by $\widetilde{\mathcal{D}}$.
     In general, $\widetilde{\mathcal{D}}$ is called the twisted Heisenberg-Virasoro algebra, and has close relationships with moduli spaces of curves\cite{ADKP}.
     The Lie algebra $\widetilde{\mathcal{D}}$ admits a Moody-Pianzola type triangular decomposition, which enables us to define its weight modules.
     Highest weight modules of $\widetilde{\mathcal{D}}$ have been well-studied in \cite{ADKP,B}.
     Moreover, simple weight modules with finite-dimensional weight spaces over $\widetilde{\mathcal{D}}$ were classified in \cite{LZ}.
     There is also a study on the support of a weight module in \cite{SS}.
     In addition to weight modules, many techniques have been developed to construct families of simple modules over $\widetilde{\mathcal{D}}$.
     In \cite{CG1}, the authors classified certain $\widetilde{\mathcal{D}}$-modules that are free of rank one over the Cartan subalgebra. 
     Untilizing simple modules of the Weyl algebra, the authors of \cite{LZG} construct new simple modules of $\widetilde{\mathcal{D}}$.
     The notion of polynomial subalgebra was introduced in \cite{M}, where it was used to construct simple modules over $\widetilde{\mathcal{D}}$. 
     Simple smooth modules of $\widetilde{\mathcal{D}}$ with non-zero level have been classified in \cite{TYZ}.

     In \cite{LP}, the infinitesimal deformations of the Lie algebra $\mathcal{D}$ were determined, resulting in a large family of deformed Lie algebras.
     The automorphism group, derivation algebra and universal central extensions of a specific deformed Lie algebra were calculated in \cite{GLP}.
     Our main interest in this paper is to study the representation theory of the universal central extension of one specific deformed Lie algebra, which shows significant differences from the original Lie algebra $\mathcal{D}$ or the universal central extension $\widetilde{\mathcal{D}}$.
     Now, let us introduce this Lie algebra.
     Let $\mathcal{L}$ be a Lie algebra with a $\C$-basis $\left\{d_{m},  h_{n}, c \mid m,n\in\mathbb{Z}\right\}$ subject to the Lie brackets
    \begin{align*}
        [d_{m}, d_{n}]&=(m-n) d_{m+n}+(m-n)h_{m+n}+\frac{m^3-m}{12}\delta _{m+n, 0}c,\\
        [d_{m}, h_{n}]&=-nh_{m+n},\quad\quad
        [h_{m}, h_{n}]=0,\quad\quad
        [\mathcal{L}, c]=0 .
    \end{align*} 
    Unlike the case of $\widetilde{\mathcal{D}}$, we show that the Lie algebra $\mathcal{L}$ does not possess a triangular decomposition in the sense of Moody and Pianzola\cite{MP}, causing difficulties in defining its weight modules.
    However, we notice that $\mathcal{L}$ remains a $\Z$-graded Lie algebra.
    A natural thought is to study its representation theory using techniques developed for $\Z$-graded Lie algebras.

    Our approach is fundamentally grounded in Mazorchuk and Zhao's classification\cite{MZ} 
    of simple Virasoro modules with local finiteness over a positive part, where their 
    explicit module constructions play a pivotal role. These foundational results not only enabled the classification of simple restricted super-Virasoro 
    modules\cite{LPX,CYZ}, but also inspired analogous constructions for extended families of algebras, including the affine Kac-Moody algebra\cite{NXZ,Ma}, planar Galilean conformal algebra\cite{GG,CY} and so on\cite{CG,MW,CHS,OW}.
    In this paper, a class of induced $\mathcal{L}$-modules will be introduced and we prove they are simple under certain conditions. 
    The final result establishes that every simple restricted $\mathcal{L}$-module of non-zero level is isomorphic to such an induced module.

    \section{Preliminaries}
    Throughout this paper, we use $\C$, $\C^*$, $\Z$, $\Z^*$, $\Z_+$ and $\N$ to denote the sets of complex numbers, non-zero complex numbers, integers, non-zero integers, non-negative integers and positive integers, respectively.
    For an arbitrary Lie algebra $\mathfrak{g}$, let $\mathcal{U}(\mathfrak{g})$ be its universal enveloping algebra.
    For $x\in\mathfrak{g}$, we use $\mathrm{ad}x:\mathfrak{g}\to\mathfrak{g}$ to denote the map that sends each $y\in\mathfrak{g}$ to $[x,y]$.

    First, we recall the definition of the $\Z$-graded Lie algebra.
    \begin{defi}
        A $\Z$-graded Lie algebra is a Lie algebra $\mathfrak{g}$ with a decomposition of subspaces
        \[
        \mathfrak{g}=\bigoplus_{n\in\Z}\mathfrak{g}_n
        \]
        such that $[\mathfrak{g}_m,\mathfrak{g}_n]\subseteq\mathfrak{g}_{m+n}$ for all $m,n\in\Z$.
        \end{defi}
        For a $\Z$-graded Lie algebra $\mathfrak{g}$, there is a notion of \textit{restricted module}.
    \begin{defi}\label{resmod}
        Let $\mathfrak{g}$ be a $\Z$-graded Lie algebra.
        A $\mathfrak{g}$-module $M$ is called a restricted module if for each $m\in M$, there exists an integer $n\in\N$ such that $\mathfrak{g}_im=0$ for all $i\ge n$.
    \end{defi} 

    Suppose $\mathfrak{g}$ is a Lie algebra.
    A linear map $\phi :\mathfrak{g}\to\mathfrak{g}$ is called an \textit{anti-automorphism} if  
    $ \phi([x, y])=[\phi(y), \phi(x)]$ for all $x, y \in \mathfrak{g}$.
    If, in addition, $\phi^2 = \mathrm{id}_{\mathfrak{g}}$,
    then $\phi$ is called an \textit{anti-involution}. 
    With the notion of anti-involution in mind,
    we now introduce the concept of \textit{Lie algebra with triangular decomposition}\cite{MP}.
    \begin{defi}\label{tri}
        Let $\mathfrak{g}$ be a Lie algebra.
        Suppose $\mathfrak{g}_{\pm}$ and $\mathfrak{h}$ are non-zero subalgebras of $\mathfrak{g}$ such that $\mathfrak{h}$ is abelian and
        \begin{equation}\label{tridecom}
            \mathfrak{g}=\mathfrak{g}_- \oplus \mathfrak{h} \oplus \mathfrak{g}_+.
        \end{equation}
        The decomposition \eqref{tridecom} is called a triangular decomposition of Moody-Pianzola type, provided that the following axioms are satisfied.
        
        \begin{enumerate}[label=(TD\arabic*)]
                \item 
                There exists an additive semigroup $Q_+ \subset \mathfrak{h}^*$ such that
                \[
                \mathfrak{g}_+=\bigoplus_{\a\in Q_+}\mathfrak{g}_{\a},
                \]
                where $\mathfrak{g}_{\a}=\left\{v\in\mathfrak{g}\mid [h,v]=\a(h)v,\forall h\in\mathfrak{h}\right\}$;
                \item\label{TD2}
                There exists an anti-involution $\sigma:\mathfrak{g}\to\mathfrak{g}$ such that $\sigma(\mathfrak{g}_+)=\mathfrak{g}_-$ and $\sigma|_{\mathfrak{h}}=\mathrm{id}_{\mathfrak{h}}$;
                \item 
                There exists a basis $\left\{\alpha_{j}\mid j \in \mathbb{J}\right\}$ of $Q_+$ consisting of $\C$-linearly independent elements of $\mathfrak{h}^*$.
                
            \end{enumerate}
    \end{defi}
    In the rest of this paper, we use the 5-tuple $(\mathfrak{g}, \mathfrak{h}, \mathfrak{g}_+, Q_+, \sigma)$ to represent a triangular decomposition of $\mathfrak{g}$.
    
    Some notions from the representation theory of Lie algebras are also needed.
    Let $\mathfrak{g}$ be a Lie algebra and $x\in\mathfrak{g}$. 
    Suppose $M$ is a $\mathfrak{g}$-module.
    The action of $x$ on $M$ is called \textit{locally nilpotent}, if for each $m\in M$, there exists some $n\in\N$ such that $x^nm=0$.
    The action of $x$ on $M$ is called \textit{locally finite}, if the dimension of the vector space $\sum_{n\in{\Z_+}}\C x^n m$ is finite.
    
    Now, we introduce the infinite-dimensional Lie algebra studied in this paper.
    
    Let $\mathcal{L}$ be a Lie algebra with the $\C$-basis $\left\{d_m,h_n,c\mid m,n\in\Z\right\}$, where the Lie brackets are as follows:
      \begin{align*}
        [d_{m}, d_{n}]&=(m-n) d_{m+n}+(m-n)h_{m+n}+\frac{m^3-m}{12}\delta _{m+n, 0}c,\\
        [d_{m}, h_{n}]&=-nh_{m+n},\quad\quad
        [h_{m}, h_{n}]=0,\quad\quad
        [\mathcal{L}, c]=0 .
    \end{align*}  
    An easy calculation shows that the center of $\mathcal{L}$ is $\C h_0 \oplus \C c$.
    We say that an $\mathcal{L}$-module has a non-zero level if $h_0$ acts on this module as a non-zero scalar.

    For $n\in\Z$, set $\mathcal{L}_n=\C d_n\oplus\C h_n\oplus\delta_{n0}\C c$.
    Then we have a decomposition
    \[
    \mathcal{L}=\bigoplus_{n\in\Z}\mathcal{L}_n
    \]
    such that $[\mathcal{L}_m,\mathcal{L}_n]\subseteq\mathcal{L}_{m+n}$ for all $m,n\in\Z$.
    It follows that $\mathcal{L}$ is a $\Z$-graded Lie algebra.
    Consequently, following Definition \ref{resmod}, an $\mathcal{L}$-module $M$ is called a restricted module if for each $m\in M$, there exists an integer $n$ (depending on $m$) such that $d_im=h_im=0$ for all $i\ge n$.

\section{Main Results}
\subsection{Triangular Decomposition}
The aim of this subsection is to prove the Lie algebra $\mathcal{L}$ does not admit a triangular decomposition.
Let $\mathcal{W}$ be the Witt algebra, that is, a Lie algebra with the basis $\left \{ L_m \mid m \in \Z\right \}$ and the Lie bracket
\[ 
[L_m, L_n]=(m-n)L_{m+n}.
\]
It is easy to see that the linear map $\tau:\mathcal{W}\to\mathcal{W}$, which maps $L_n$ to $L_{-n}$, is an anti-involution of $\mathcal{W}$.
According to \cite{DZ},
if $\phi:\mathcal{W} \to \mathcal{W}$ is an automorphism of Lie algebra,
then $\phi(L_n)=\varepsilon \lambda^nL_{\varepsilon n}$ for some $\lambda \in \C^*$ and $\varepsilon \in\left \{ \pm1 \right \} $.
Thus, each anti-automorphism of $\mathcal{W}$ is of the form
\begin{equation}\label{anti-auto}
\begin{split}
\tau\phi:\mathcal{W}&\to\mathcal{W},\\
L_n&\mapsto \epsilon\lambda^nL_{-\epsilon n}.
\end{split}
\end{equation}
The following lemma is due to Chari and Pressley\cite{CP}.
\begin{lemma}\label{W_0}
Suppose $0 \ne x \in \mathcal{W}$ such that $\mathrm{ad}x$ is diagonalizable on $\mathcal{W}$, that is,
\[\mathcal{W}=\bigoplus_{\mu \in \C} \mathcal{W}_{\mu}^{x}, \]
where for $\mu\in\C$, $\mathcal{W}_{\mu}^{x}:=\left \{w \in  \mathcal{W} \mid [x, w]=\mu w\right \}$.
Then there exists some $\xi \in \C^*$ such that $x=\xi L_0$.
\end{lemma}
This lemma allows us to determine all the triangular decompositions of $\mathcal{W}$.
\begin{prop}\label{tri of w}
   Suppose $\left(\mathcal{W}, \mathfrak{H}, \mathcal{W}_+, Q_+, \sigma\right)$ is a triangular decomposition of $\mathcal{W}$.
    Then we have
    \[
    \mathfrak{H}=\C L_0,\quad 
    \mathcal{W}_+=\bigoplus_{n \in \N}\C L_{-n},\quad 
     Q_+=\N,\quad\mbox{and}\quad\sigma(L_n)=\lambda^{n}L_{-n}\quad\mbox{for some }\lambda\in\C^*.
    \]
\begin{proof}
From the decomposition \eqref{tridecom}, the adjoint action of every $0\ne x \in \mathfrak{H}$ is semisimple.
Hence, we have $\mathfrak{H}=\C L_0$ by Lemma $\ref{W_0}$.
Assume that $\sum_{i=p}^{q} \lambda_iL_i\in\mathcal{W}$
is an eigenvector of $\mathrm{ad}L_0$ with respect to the eigenvalue $\mu\in\C$, where $\lambda_q \neq 0$.
Then we have 
\[ 
[L_0, \sum _{i=p}^{q}\lambda_i L_i]= \mu\sum _{i=p}^{q}\lambda_i L_i,
\]
that is,
\[ 
\sum _{i=p}^{q}(\mu+i)\lambda_i L_i=0.
\]
It follows that $\lambda_i = 0$ for all $i \neq q$,
which implies $\mathcal{W}_+=\bigoplus_{n \in \N}\C L_{-n}$ and $Q_+\cong\N$.
Combining \eqref{anti-auto} with the fact that $\sigma(\mathcal{W}_+)=\mathcal{W}_-$,
we obtain $\sigma(L_n)=\lambda^{n}L_{-n}$ for some $\lambda\in\C^*$.
Now, it is easy to see that the 5-tuple $(\mathcal{W}, \mathfrak{H}, \mathcal{W}_+, Q_+, \sigma)$ is a triangular decomposition of $\mathcal{W}$.
\end{proof}
\end{prop}
Next, we turn to study the Lie algebra $\mathcal{L}$.
Let $\mathcal{I}$ be the subspace of $\mathcal{L}$ spanned by $\left\{ h_m, c\mid m\in\Z\right\}$.
Then, it is straightforward to check $\mathcal{I}$ is the unique maximal ideal of $\mathcal{L}$.
It follows that $\tau(\mathcal{I})=\mathcal{I}$ for each anti-involution $\tau$.
\begin{prop}
    The Lie algebra $\mathcal{L}$ does not admit a Moody-Pianzola type triangular decomposition.
    \begin{proof}
        Suppose, for contradiction, that $(\mathcal{L}, \mathcal{H}, \mathcal{L}_+, Q_+, \tau)$ is a triangular decomposition.
        Then we have the decomposition 
        \[ 
        \mathcal{L} /\mathcal{I} = (\mathcal{L}_-+\mathcal{I})/\mathcal{I} \ \oplus  ( \mathcal{H}+\mathcal{I})/\mathcal{I} \oplus( \mathcal{L}_++\mathcal{I})/\mathcal{I}.
        \]
        If $\mathcal{L}_+  \subseteq \mathcal{I}$, then the inclusion $\mathcal{L}_-=\tau(\mathcal{L}_+)\subseteq \tau(\mathcal{I})=\mathcal{I}$ holds.
        Hence the quotient algebra $\mathcal{L}/\mathcal{I}=(\mathcal{H}+\mathcal{I})/\mathcal{I}$ is abelian, which contradicts to the fact that $\mathcal{L}/\mathcal{I} \cong \mathcal{W}$.
        If $\mathcal{H} \subseteq \mathcal{I}$,
        then we have $\mathcal{L}_\pm=[\mathcal{H},\mathcal{L}_\pm]\subseteq [\mathcal{I},\mathcal{L}]\subseteq \mathcal{I}$, and thus $\mathcal{L}=\mathcal{I}$, which is a contradiction. 
        Since the adjoint action of any non-zero element in $(\mathcal{H}+\mathcal{I})/\mathcal{I}$ on $\mathcal{L}/\mathcal{I}$ is semisimple,
        using Lemma \ref{W_0} and the Lie algebra isomorphism $\mathcal{L}/\mathcal{I} \cong \mathcal{W}$,
        we have 
        \begin{equation} \label{H}
        (\mathcal{H}+\mathcal{I})/\mathcal{I}=\C(d_0+\mathcal{I}).
        \end{equation}

        Since $\mathcal{L}_+ \not\subseteq \mathcal{I}$,
        there exist $\alpha \in Q_+$, $x\in\mathcal{H}$, $y\in\mathcal{L}_{\alpha}\setminus \mathcal{I}$ and $\zeta\in\C^*$ such that 
        \begin{equation}\label{x,y}
         [x, y]=\zeta y.   
        \end{equation}
        Passing to the quotient algebra $\mathcal{L}/\mathcal{I}$,
        we have
        $
        [x+\mathcal{I}, y+\mathcal{I}]=\zeta (y+\mathcal{I}),
        $
        which shows that $x \notin  \mathcal{I}$.
        Using the Lie algebra isomorphism $\mathcal{L}/\mathcal{I} \cong \mathcal{W}$, 
        we know that
        \[
        x-\alpha d_0,\; y-\beta d_q\in\mathcal{I}
        \]
        for some $q\in\Z^*$ and $\alpha,\beta\in\C^*$.
        It follows that 
        \[
        x=\alpha d_0+\sum_{i\in\Z}\alpha_ih_i+\mu c \quad\quad \text{and}\quad\quad
        y=\beta d_q+\sum_{j\in\Z}\beta_{j}h_j +\eta c,
        \]
        where $\mu,\eta\in\C$ and only finitely many $\a_i$ and $\beta_j$ are non-zero.
        According to equation \eqref{x,y}, we have 
        \[
       (\zeta+q\alpha ) \beta d_q+\sum_{j \in \Z}[(\zeta+j \alpha )\beta_j-(j-q)\beta\alpha_{j-q}+\delta_{jq}q\alpha \beta]h_j +\zeta \eta c=0.
        \]
        Extracting the coefficients of $d_q$ and $h_q$,
        we obtain $(\zeta+q\alpha  )\beta=0 $ and $(\zeta+q\alpha)\beta _q+q\alpha \beta =0$, respectively.
        This implies $q \alpha \beta=0$, contradicting to the fact that $q\in\Z^*$ and $\alpha,\beta\in\C^*$.
        Hence $\mathcal{L}$ does not admit a triangular decomposition.   
    \end{proof}
\end{prop}

\subsection{Simple Induced Modules}
    In this subsection, we introduce a class of induced modules over $\mathcal{L}$ and prove they are simple under certain conditions.
    These simple modules will be useful in the classification of simple restricted $\mathcal{L}$-modules.

    For $q\in \mathbb{Z}_+\cup\left\{+\infty\right\}$, define the subalgebra $\mathcal{L}_q$ of $\mathcal{L}$ as follows:
    \begin{equation}
    \label{subalg}
    \mathcal{L}_{q}:=\sum_{i \in \mathbb{Z_{+}}}(\mathbb{C} d_{i}+\mathbb{C} h_{i-q})+\C c.
    \end{equation}
    For an arbitrary $\mathcal{L}_{q}$-module $V$, we have an induced $\mathcal{L}$-module 
    $$\operatorname{Ind}_q(V):=\mathcal{U}(\mathcal{L}) \otimes_{\mathcal{U}\left(\mathcal{L}_{q}\right)}V.$$

    Let $\mathbb{M}$ be the set of infinite sequences of the form $\boldsymbol{i} =(\ldots, i_{3}, i_{2},i_{1})$, where each $i_k\in\Z_+$ and only finitely many of them are non-zero.
    Let $\boldsymbol{0}$ be the sequence with all entries equal to zero.
    For $i\in\N$, denote by $\boldsymbol{\epsilon}_i$ the sequence $( \ldots,\delta_{i3},\delta_{i2},\delta_{i1})$.
    For any $\boldsymbol{i} \in \mathbb{M}$, define
    $$ 
    \operatorname{w}(\boldsymbol{i}):=\sum_{s \in \mathbb{Z}_{+}} s \cdot i_{s}.
    $$
    For any $\boldsymbol{0} \ne \boldsymbol{i} \in \mathbb{M}$, let $p$ (respectively $q$) be the largest (respectively smallest) integer such that $i_{p} \neq 0$ (respectively $i_{q} \neq 0$).
    Denote by $\boldsymbol{i}'= \boldsymbol{i}-\boldsymbol{\epsilon}_{p}$ and $\boldsymbol{i}''= \boldsymbol{i}-\boldsymbol{\epsilon}_{q}$, respectively.

        For $\boldsymbol{i}, \boldsymbol{j} \in \mathbb{M}$, define the \textit{lexicographical total order} $>$ and the \textit{reverse lexicographical total order} $\succ $ on $\mathbb{M}$ as follows.
        \begin{align*}
        \boldsymbol{j}>\boldsymbol{i} &\Leftrightarrow\text{There exists } r \in \mathbb{N} \text { such that }j_{s}=i_{s}\mbox{ for all } s>r\text{ and }j_{r}>i_{r}.\\
        \boldsymbol{j} \succ \boldsymbol{i} &\Leftrightarrow \mbox{There exists } r \in \mathbb{N} \text { such that }j_{s}=i_{s}\mbox{ for all } 1 \leqslant s<r\text { and } j_{r}>i_{r}.
        \end{align*}
        Furthermore,
        there is an induced total order on $\mathbb{M} \times \mathbb{M}$, still denoted by $ \succ$, as follows:
        \begin{equation}
        \label{order}
        (\boldsymbol{i}, \boldsymbol{j}) \succ (\boldsymbol{k}, \boldsymbol{l})  \Leftrightarrow (\boldsymbol{j}, \operatorname{w}(\boldsymbol{j})) \succ(\boldsymbol{l}, \operatorname{w}(\boldsymbol{l})) \text { or }(\boldsymbol{j}=\boldsymbol{l} \text { and } \boldsymbol{i}>\boldsymbol{k}),\quad  
        \mbox{ where } \boldsymbol{i}, \boldsymbol{j}, \boldsymbol{k}, \boldsymbol{l} \in \mathbb{M}.
        \end{equation}

        For $\boldsymbol{i}, \boldsymbol{j} \in \mathbb{M}$, denote by
        $$ h^{\boldsymbol{i}} d^{\boldsymbol{j}}=\cdots h_{-q-2}^{i_{2}} h_{-q-1}^{i_{1}} \cdots d_{-2}^{j_{2}} d_{-1}^{j_{1}} \in \mathcal{U}(\mathcal{L}).$$
        According to the Poincar\'{e}-Birkhoff-Witt Theorem,
        every $v\in\operatorname{Ind}_q(V)$ can be uniquely written as
        \begin{equation}\label{o}
        \sum_{\boldsymbol{i}, \boldsymbol{j} \in \mathbb{M}} h^{\boldsymbol{i}} d^{\boldsymbol{j}} v_{\boldsymbol{i}\boldsymbol{j}},    
        \end{equation}
        where $v_{\boldsymbol{i}\boldsymbol{j}} \in V $ and only finitely many of them are non-zero.
        For every $0\ne v\in\operatorname{Ind}(V)$ written in the form of \eqref{o}, define the \textit{support} of $v$ to be
        $$\operatorname{supp}(v):=  \left \{ (\boldsymbol{i}, \boldsymbol{j}) \in \mathbb{M} \times \mathbb{M} \mid  v_{\boldsymbol{i}\boldsymbol{j}} \neq 0  \right \}.$$
        Define the \textit{degree} of $v$ to be the maximal element of $\operatorname{supp}(v)$ with respect to the total order on $\mathbb{M} \times \mathbb{M}$ defined in \eqref{order}, and denote it by $\deg(v)$.

        Now, we introduce the following technical lemma.
    \begin{lemma}
        \label{b}
            Let $q \in \mathbb{Z}_+ \cup\{+\infty\}$ and $V$ be an $\mathcal{L}_{q}$-module (not necessarily simple).
            Assume that there exists some $k\in\Z_+$ such that $h_{k}$ acts injectively  on $V$ and $ h_{i}V=d_{j}V=0$ for all $ i > k $ and $ j > k + q $. 
            For any $v \in \operatorname{Ind}(V) \setminus V$,
            suppose $\deg(v)=(\boldsymbol{i}, \boldsymbol{j})$.
            Let $a=\max \left\{s: i_{s} \neq 0\right\}$ if $ \boldsymbol{i}  \neq 0$ and $b=\min \left\{s: j_{s} \neq 0\right\} $ if $ \boldsymbol{j}  \ne 0$.
            Then the following statements hold.
            \begin{enumerate}[label=(\arabic*)]
                \item 
                If $ \boldsymbol{j} \ne \mathbf{0}$, then $b > 0$ and $\operatorname{deg}(h_{b+k}v)=(\boldsymbol{i}, \boldsymbol{j}'')$.
                \item 
                If $ \boldsymbol{j} = \boldsymbol{0}$ and $ \boldsymbol{i} \neq \boldsymbol{0}$, then $a > 0$ and $\operatorname{deg}(d_{a+k+q}v)=(\boldsymbol{i}', 0)$.
            \end{enumerate}
            \begin{proof}
            (1) For any $(\boldsymbol{x}, \boldsymbol{y}) \in \operatorname{supp}(v)$,
            using the Lie brackets 
            \[
            [d_{m}, h_{n}]=-nh_{m+n}\quad\quad \mbox{and}\quad\quad[h_{m}, h_{n}]=0,
            \]
            we have 
            \begin{equation}\label{a}
            \begin{aligned}
            h_{b+k}h^{\boldsymbol{x}}d^{\boldsymbol{y}}v_{\boldsymbol{x}\boldsymbol{y}} & =h^{\boldsymbol{x}}[h_{b+k}, d^{\boldsymbol{y}}]v_{\boldsymbol{x} \boldsymbol{y}} \\
            & =\sum_{\boldsymbol{z} \in \mathbb{I} }\lambda_{\boldsymbol{z}}h^{\boldsymbol{x}}d^{\boldsymbol{z}}h_{\varphi (\boldsymbol{z})}v_{\boldsymbol{x}\boldsymbol{y}},
            \end{aligned}
            \end{equation}  
            where $\lambda_{\boldsymbol{z}}\in\C^*$ and $$\mathbb{I}=\left \{ \boldsymbol{z} \in \mathbb{M} \mid z_{i} \le y_{i}  \mbox{ for each }  i \in \N,\mbox{ }\varphi (\boldsymbol{z})\in\Z_{\le k},\mbox{ and }\operatorname{w}(\boldsymbol{z})-\operatorname{w}(\boldsymbol{y})=\varphi(\boldsymbol{z})-b-k \right \}.$$

            First, suppose $\operatorname{w}(\boldsymbol{y}) = \operatorname{w}(\boldsymbol{j})$.
            If $\boldsymbol{y}=\boldsymbol{j}$,
            then according to equation \eqref{a},
            we have 
                \begin{align*}
                h_{b+k}h^{\boldsymbol{x}}d^{\boldsymbol{j}}v_{\boldsymbol{x}\boldsymbol{j}} & 
                 = h^{\boldsymbol{x}}[h_{b+k}, d^{\boldsymbol{j}}]v_{\boldsymbol{x}\boldsymbol{j}} \\
                 &=j_{b}(b+k)\cdot h^{\boldsymbol{x}}d^{\boldsymbol{j}''}(h_{k}  v_{\boldsymbol{x}\boldsymbol{j}})+\sum_{\boldsymbol{z} \in \mathbb{H} }\lambda_{ \boldsymbol{z}} h^{\boldsymbol{x}}d^{\boldsymbol{z}}h_{\varphi (\boldsymbol{z})}v_{\boldsymbol{x} \boldsymbol{j}},
                \end{align*}
            where $$\mathbb{H}=\left \{ \boldsymbol{z} \in \mathbb{M} \mid z_{i} \le j_{i}  \mbox{ for each }  i \in \N\mbox{ and } \operatorname{w}(\boldsymbol{z}) < \operatorname{w}(\boldsymbol{j})-b\right \}.$$
            It follows that  $ \operatorname{deg}(h_{b+k}h^{\boldsymbol{x}}d^{\boldsymbol{j}}v_{\boldsymbol{x} \boldsymbol{j}}) \preceq  (\boldsymbol{i}, \boldsymbol{j}'')$, where the equality holds if and only if 
            $\boldsymbol{x}=\boldsymbol{i}$.
            If $\boldsymbol{y} \prec \boldsymbol{j}$, then
            according to equation \eqref{a},
            we have 
            $$\operatorname{w}(\boldsymbol{z})=\operatorname{w}(\boldsymbol{y})-b+\varphi(\boldsymbol{z})-k\le\operatorname{w}(\boldsymbol{y})-b = \operatorname{w}(\boldsymbol{j})-b =\operatorname{w}(\boldsymbol{j}'').$$
            If $  \operatorname{w}(\boldsymbol{z}) < \operatorname{w}(\boldsymbol{y})-b $, then the proof is trivial.
            Assume that $\operatorname{w}(\boldsymbol{z}) = \operatorname{w}(\boldsymbol{y})-b$, that is, $$\sum_{i\in\N}iz_i=\sum_{i\in\N}iy_i-b.$$
            If $y_b=0$, then using the fact $\boldsymbol{y} \prec \boldsymbol{j}$ and $z_i\le y_i$ for each $i\in\N$, 
            we obtain $b=\sum_{i\ge b+1}i(y_i-z_i)\ge b+1$, which is a contradiction.
            Hence, 
            we have $\boldsymbol{z}=\boldsymbol{y}'' \prec \boldsymbol{j}''$.

            Finally, suppose $w(\boldsymbol{y}) < w(\boldsymbol{j})$.
            From equation \eqref{a},
            we have $\operatorname{w}(\boldsymbol{z})-\operatorname{w}(\boldsymbol{y})=\varphi(\boldsymbol{z})-b-k$.
            It follows that   $\operatorname{w}(\boldsymbol{z}) \le  \operatorname{w}(\boldsymbol{y}) - b< \operatorname{w}(\boldsymbol{j})-b = \operatorname{w}(\boldsymbol{j}'')$, which implies $\operatorname{deg}(h_{b+k}h^{\boldsymbol{x}}d^{\boldsymbol{y}}v_{\boldsymbol{x}\boldsymbol{y}}) \prec (\boldsymbol{i}, \boldsymbol{j}'') $.


            (2) If  $ \boldsymbol{j} = 0$ and $ \boldsymbol{i} \neq 0$,
            we deduce that
            \begin{center}
            \begin{math}
                \begin{aligned}
                d_{a+k+q} h^{\boldsymbol{i}}v_{\boldsymbol{i}\boldsymbol{0}} & = [ d_{a+k+q}, h_{-q-a}^{i_{a}} \cdots h_{-q-1}^{i_{1}}  ] v_{\boldsymbol{i}\boldsymbol{0}} \\
                &=i_{a}(q+a)\cdot h_{-q-a}^{i_{a}-1} \cdots h_{-q-1}^{i_{1}}h_kv_{\boldsymbol{i}\boldsymbol{0}}. \\
               \end{aligned}
            \end{math}
            \end{center} 
            It follows that $\operatorname{deg}(d_{a+k+q} h^{\boldsymbol{i}}v_{\boldsymbol{i}\boldsymbol{0}})=(\boldsymbol{i}', 0)$.
            For $(\boldsymbol{x}, 0) \prec (\boldsymbol{i}, 0)$,
            if $ x_{a}=0$,
            then, by direct computations one has $d_{a+k+q}h^{\boldsymbol{x}}v_{\boldsymbol{x}\boldsymbol{0}}=0$,
            otherwise, $ 0 \neq x_{a}< i_{a}$, then it is clear that $\operatorname{deg}(d_{a+k+q}h^{\boldsymbol{x}}v_{\boldsymbol{x}\boldsymbol{0}})=(\boldsymbol{x}', 0)\prec (\boldsymbol{i}', 0) $. 
            \end{proof}
    \end{lemma}
    \begin{rmk}
    Assume that $U$ is a simple $\mathcal{L}_q$-module and $ h_{i}U=0$ for all $ i> k$. 
    Then, either $h_kU=0$ or the action of $h_k$ on $U$ is injective.
    This is because $\mathrm{Ker}(h_k)$ is an $\mathcal{L}_q$-submodule of $U$. 
    \end{rmk}
    Let $V$ be a simple $\mathcal{L}_q$-module satisfying the hypothesis of Lemma \ref{b}. 
    For an arbitrary element $0\ne u\in \operatorname{Ind}_q(V)\setminus V$.
    Utilizing Lemma \ref{b} repeatedly, we obtain a non-zero element lying in $V$.
    Hence, we obtain the \textbf{main theorem} of this subsection.
    \begin{thm}\label{t1}
    Suppose $q \in \mathbb{Z}_+  \cup \left \{ +\infty  \right \} $ and $V$ is a simple $\mathcal{L}_q$-module. 
    Assume that there exists $ k \in \mathbb{Z}_{+}$ such that
         $h_{k}V\ne 0$ and
         $d_{m}V=h_{n}V=0$ for all $m>k+q$ and $n>k$.
    Then the induced $\mathcal{L}$-module $\mathrm{Ind}_q(V)$ is simple.
    \end{thm}
    \begin{rmk}
    Let $\mathcal{L}_0$ be the subalgebra of $\mathcal{L}$, as defined in \eqref{subalg}.
    Define a one-dimensional module $\C v_{\lambda\mu}$ of $\mathcal{L}_0$ with $c$ and $h_0$ acting as scalars $\lambda,\mu\in\C$, respectively.
    The induced $\mathcal{L}$-module $\mathrm{Ind}_0(\C v_{\lambda\mu})$ is called the \textit{generalized Verma module} of $\mathcal{L}$.
    If $\mu\ne0$, then from Theorem \ref{t1}, the $\mathcal{L}$-module $\mathrm{Ind}_0(\C v_{\lambda\mu})$ is simple. Conversely, suppose $\mu=0$, then by direct computations one has
    $$
    d_m(h_{-1}\otimes v_{\lambda\mu}),h_m(h_{-1}\otimes v_{\lambda\mu})\in\C(h_{-1}\otimes v_{\lambda\mu})
    $$
    for all $m\in\Z_+$.
    It follows that $h_{-1}\otimes v_{\lambda\mu}\in\mathrm{Ind}_0(\C v_{\lambda\mu})$ generates a non-zero proper submodule of $\mathrm{Ind}_0(\C v_{\lambda\mu})$.
    Therefore, the generalized Verma module $\mathrm{Ind}_0(\C v_{\lambda\mu})$ is simple if and only if $\mu\ne0$. 
    \end{rmk}
        
    \subsection{Simple Restricted Modules of Non-zero Level}
    In this subsection, we shall classify simple restricted $\mathcal{L}$-modules of non-zero level.
    The following elementary lemma will be useful.
    \begin{lemma}\label{l3}
     Let $\mathfrak{g}$ be a Lie algebra with a basis $\left \{ x_1, x_2, \ldots\right \}$ and $V$ be a $\mathfrak{g}$-module.
     Suppose there exists a non-zero vector $v \in V$ such that $\sum_{i=1}^{+\infty} \C x_i v=\sum_{i=1}^{m} \C x_i v$.
     Then the space $W=\sum_{n_1, n_2, \cdots, n_m \in \Z_+}\C x_{1}^{n_1}x_{2}^{n_2}\cdots x_{m}^{n_m}v$ is a $\mathfrak{g}$-submodule of $V$.
    \end{lemma}
    Using the above lemma, we now characterize the simple induced $\mathcal{L}$-module $\operatorname{Ind}_q(V)$.
    \begin{thm}\label{d}
            Suppose $S$ is a simple $\mathcal{L}$-module with $h_0$ acting as a non-zero scalar.
            Then the following conditions are equivalent:
            \begin{enumerate}[label=(\arabic*)]
            \item 
            There exists $t \in \mathbb{Z}_+ $ such that  the actions of $d_{i}$, $h_{i}$, $i \geq t $ on $S$ are locally finite.
            \item 
            There exists $t \in \mathbb{Z}_+ $ such that  the actions of $d_{i}$, $h_{i}$, $i \geq t $ on $S$ are locally nilpotent.
            \item 
            There exists $q \in \mathbb{Z}_+$ and a simple $\mathcal{L}_{q}$-module $V$ satisfying the conditions in Theorem \ref{t1} such that $ S \cong \operatorname{Ind}_q(V)$. 
            \end{enumerate}

            \begin{proof}
            $(3)\Rightarrow (2)\Rightarrow (1)$ is obvious.
            Now, we prove $(1)\Rightarrow (3)$.
            
            For $s\in\Z_+$, let $\mathcal{L}^s$ be the subalgebra of $\mathcal{L}$ generated by $\left\{d_m,h_m\mid m\in\Z_{\ge s}\right\}$.
            Suppose $S$ is a simple $\mathcal{L}$-module.
            Assume that there exists some $t \in \N$ such that the actions of $d_{i}$ and $h_{i}$ on $S$ are locally finite for all $i\ge t$.
            Then, we can find $0\ne v \in S$ such that $d_{t}v=\lambda v$ for some $\lambda \in \C$.

            For any integer $j > t $, denote by 
            \[ 
            N_{dj}=\sum_{n \in \Z_{+}} \C d_{t}^{n} d_{j}v = \mathcal{U}(\C d_{t})d_{j}v \quad\mbox{and}\quad   N_{hj}=\sum_{n \in \Z_{+}} \C d_{t}^{n} h_{j}v = \mathcal{U}(\C d_{t})h_{j}v, \quad\mbox{respectively}.
            \]
            Owing to (1), both $N_{dj}$ and $N_{hj}$ are finite-dimensional.
            
            We claim that $h_{j+mt}v \in N_{hj}$ and $d_{j+mt}v \in N_{dj}+N_{hj}$ for all $m \in \Z_{+}$.
            In the case of $m=0$, this certainly holds true.
            Assume that $h_{j+at}v \in N_{hj}$  and $d_{j+at}v \in N_{dj}+N_{hj}$ for all $a \le m$.
            Then, we have 
                \begin{align*}
                    -(j+mt)h_{j+(m+1)t}v& =[d_{t},h_{j+mt}]v \\
                                            & =d_{t}h_{j+mt}v-h_{j+mt}d_{t}v \\
                                            & =(d_{t}-\lambda)h_{j+mt}v
                \end{align*}
                and
                \begin{align*}
                    (t-j-mt)d_{j+(m+1)t}v+(t-j-mt)h_{j+(m+1)t}v& =[d_{t},d_{j+mt}]v \\
                                            & =d_{t}d_{j+mt}v-d_{j+mt}d_{t}v \\
                                            & =(d_{t}-\lambda)d_{j+mt}v.
                \end{align*}
                From the two equations above,
                we obtain that $h_{j+(m+1)t}v \in N_{hj}$  and $d_{j+(m+1)t}v \in N_{dj}+N_{hj} $ since $j+mt, t-j-mt \neq 0$.

            Now,
            we see that $\sum_{m \in \Z_+}\C d_{j+mt}v \subseteq N_{dj}+N_{hj}$ and $\sum_{m \in \Z_+}\C h_{j+mt}v \subseteq N_{hj}$.
            Hence,
            they are both finite-dimensional for all $j>t$. 
            It follows that 
            \[
            \sum_{i \in \Z_+}\C d_{t+i}v=\C d_{t}v +\sum_{s=t+1}^{2t}\left(\sum_{n \in \Z_+}\C d_{s+nt}v\right) 
            \]
            and
            \[\sum_{i \in \Z_+}\C h_{t+i}v=\C h_{t}v +\sum_{s=t+1}^{2t}\left(\sum_{n \in \Z_+}\C h_{s+nt}v\right) 
            \]
            are finite-dimensional.
            Therefore, we can find some $l \in \N $ such that 
            \[\sum_{i \in \Z_+}\C d_{t+i}v=\sum_{i=0}^{l}\C d_{t+i}v \quad\mbox{and}\quad \sum_{i \in \Z_+}\C h_{t+i}v=\sum_{i=0}^{l}\C h_{t+i}v.  \]
            Combining the assumptions of (1) and Lemma \ref{l3},
            we know that 
            $$W=\sum_{m_i,n_i\in\Z_+} \C d_{t}^{m_0}d_{t+1}^{m_1} \cdots d_{t+l}^{m_l}h_{t}^{n_0}h_{t+1}^{n_1} \cdots h_{t+l}^{n_l}v$$
            is a finite-dimensional $\mathcal{L}^{t}$-module.

            Let $\mathrm{Ann}(W)$ be the annihilator of the $\mathcal{L}^{t}$-module $W$.
            Since $W$ is finite-dimensional,
            we must have $\mathrm{Ann}(W) \neq 0$.
            Assume that there exists 
            \[
            0\ne w=\sum_{i=0}^{m}a_id_{t+i}+\sum_{i=0}^{m}b_ih_{t+i}\in\mathrm{Ann}(W),
            \]
            where $a_i,b_i\in\C$.
            Applying $h_t$ and note that $\mathrm{Ann}(W)$ is an ideal of $\mathcal{L}^t$, we see that
            \[
            \sum_{i=0}^{m}ta_ih_{2t+i}=[h_t,w]\in\mathrm{Ann}(W).
            \]
            Now, choose a minimal $n\in\Z_+$ such that
            \[
            w_1:=\sum_{i=0}^{n}c_ih_{r+i}\in\mathrm{Ann}(W),
            \]
            where $c_0\ne 0$ and $r\in\Z_{\ge t}$.
            Then, we have
            \[
            w_2:=\sum_{i=0}^{n}c_i(t+r+i)(r+i)h_{2t+r+i}=[[w_1,d_t],d_t]\in\mathrm{Ann}(W)
            \]
            and
            \[
            w_3:=\sum_{i=0}^{n}c_i(r+i)h_{2t+r+i}=[w_1,d_{2t}]\in\mathrm{Ann}(W).
            \]
            It follows that $$\sum_{i=1}^{n}i(r+i)c_ih_{2t+r+i}=w_2-(t+r)w_3\in\mathrm{Ann}(W).$$
            Hence, we see that $c_i=0$ for $1\le i\le n$ and $h_rW=0$.
            Now, for $m\in\Z_{\ge t}$, we have 
            $$h_{r+m}W=\frac{1}{r}[h_r,d_m]W=0.$$
            Let $\mathcal{I}^{r+t}=\sum_{i\in\Z_+}\C h_{r+t+i}$ and $\mathcal{L}^{r+t}/\mathcal{I}^{r+t}$ be the corresponding quotient Lie algebra.
            Then $W$ is a module over $\mathcal{L}^{r+t}/\mathcal{I}^{r+t}$ since $\mathcal{I}^{r+t}W=0$.
            Similar to the above discussions, we see that $d_i W=0$ when $i$ is sufficiently large.
        
            For $x, y \in \Z$, consider the vector space 
            \[ N_{xy}:=\left\{v \in S \mid h_{i} v=d_{j} v=0 \text { for all } i>x \text { and } j>y\right\}.\]
            From the discussions above, there exists a sufficiently large integer $x_0$ such that $0\ne W \subseteq N_{x_0, x_0+1}$.
            Since $h_0$ acts on $S$ as a non-zero scalar, it follows that $N_{-1, x_0+1}=0$.
            Now, we have a chain of vector spaces
            \[ 
            0=N_{-1, x_0+1}\subseteq N_{0, x_0+1}\subseteq \cdots\subseteq N_{x_0-1, x_0+1}\subseteq N_{x_0, x_0+1} \ne0.
            \]
            Hence there exists $-1 \le k \le x $ such that $N_{k-1, x_0+1}=0$ and $N_{k, x_0+1} \neq 0$, which shows that $h_{k}$ acts injectively on $N_{k, x_0+1}$.
            Let $q=x_0+1-k$ and $V=N_{k, x_0+1}$, respectively.
            It is obvious that $V$ is an $\mathcal{L}_{q}$-module.

            According to the fact that $S$ is a simple $\mathcal{L}$-module,
            there exists a canonical surjective map
            $
            \pi: \operatorname{Ind}_q(V) \to S$ such that $\pi(1 \otimes v)=v$ for all $v\in V$.
            Now we prove the map $\pi$ is injective.
            Denote by $K:=\mathrm{Ker}(\pi)$. 
            It is apparent that $K\cap V=0$.
            For $u\in K$, if $u\ne 0$, then by the same arguments in the proof of Lemma \ref{b},
            we shall find an element $0\ne u' \in K\cap V$.
            This leads to a contraction.
            Hence, we have $K=0$ and the map $\pi$ is bijective.
            By the Poincar\'{e}-Birkhoff-Witt Theorem,
            the universal enveloping algebra $\mathcal{U}(\mathcal{L})$ is flat as a right $\mathcal{U}(\mathcal{L}_q)$-module.
            Therefore, the tensor product functor $\mathcal{U}(\mathcal{L})\otimes_{\mathcal{U}(\mathcal{L}_q)}\underline{\;\;\,\;}$ is exact.
            It follows that $V$ is a simple $\mathcal{L}_q$-module.
            The proof is completed.
    \end{proof}
    \end{thm}
    The following lemma admits a straightforward proof, which we omit for brevity.
    \begin{lemma}\label{c}
        Suppose $S$ is a simple restricted $\mathcal{L}$-module.
        Then there exists $s\in\mathbb{N}$ such that the actions of $d_{t}$ and $h_{t}$ on $S$ are locally nilpotent for all $t \ge s$.
    \end{lemma}
    Combining Theorem \ref{d} and Lemma \ref{c}, we obtain the \textbf{main theorem} of this paper.
    \begin{thm}\label{main}
    Every simple restricted $\mathcal{L}$-module of non-zero level is isomorphic to $\mathrm{Ind}_q(V)$ for some $q\in\Z_+$ and simple $\mathcal{L}_q$-module $V$.
    \end{thm}

	\section*{Acknowledgements}
	This paper is partially supported by the National Key R\&D Program of China (2024YFA1013802), NSF of China (11931009, 12101152, 12161141001, 12171132, 12371040, and 12401036), Innovation Program for Quantum Science and Technology (2021ZD0302902), and the Fundamental Research Funds for the Central Universities of China.

\end{document}